\documentclass{article}
\usepackage{amsfonts}
\usepackage{amssymb}
\usepackage{amsmath}

\input{tcilatex}
\begin{document}

\author{Piotr WILCZEK}
\title{Under the Continuum Hypothesis all nonreflexive Banach space \\
ultrapowers are primary.}
\date{}
\maketitle

\begin{abstract}
In this note a large class of primary Banach spaces is characterized.
Namely, it will be demonstrated that - under the Continuum Hypothesis - if $%
E $ is any infinite dimensional nonsuperreflexive Banach space, then its
ultrapower $\ell _{\infty }(E)/\mathcal{N}_{\mathcal{U}}$ is always primary.
Consequently, any infinite dimensional nonsuperreflexive Banach space can be
isometrically embedded into its primary ultrapowers.
\end{abstract}


$2010$ \textit{Mathematical Subject Classification}. $46B08$, $46B20$, $%
46B25 $.\bigskip

$\mathbf{1.}$\textbf{\ Introduction}.\bigskip

Recall that a Banach space $E$ is called \textit{primary} if, for every
bounded projection operator $P$ on $E$, either $PE$ or $(I-P)E$ is
isomorphic to $E$ $($where $I$ is the identity operator$).$ It was proved
that many classical function spaces $C(K)$, i.e., the space of continuous
scalar-valued functions on any infinite metrizable compact space $K$ and $%
L_{p}(0,1)$ spaces $($for $1\leq p\leq \infty )$ are all primary. Also it
was demonstrated that Pe\l czy\'{n}ski's universal basis space $U_{1}$ and
the space $J$ of James are primary $($\cite{C, D, DR, LT}$)$.

In this note the methodology of Banach space ultrapowers will be employed in
order to single out a large class of Banach spaces which are primary.
Namely, it will be shown that if the C\textit{ontinuum Hypothesis} $(CH)$
holds and $E$ is any infinite dimensional nonsuperreflexive Banach space,
then the ultrapower $\ell _{\infty }(E)/\mathcal{N}_{\mathcal{U}}$ is always
primary. Therefore, any infinite dimensional nonsuperreflexive Banach space $%
E$ can be isometrically embedded into the primary Banach space of the form $%
\ell _{\infty }(E)/\mathcal{N}_{\mathcal{U}}$.\bigskip

$\mathbf{2}$\textbf{. Banach space ultrapowers and the Representation
Theorem\bigskip }

Recall that a Banach space $E$ is said to be finite dimensional if and only
if its unit ball is compact, i.e., if and only if for every bounded family $%
(x_{i})_{i\in I}$ and for every ultrafilter $\mathcal{U}$ on the set $I$,
the so-called $\mathcal{U-}$limit $\underset{i,\mathcal{U}}{\lim }x_{i}$
exists. The sequence $(x_{i})_{i\in I}$ converges to $x$ with respect to the
ultrafilter $\mathcal{U}$ if and only if the set $\left\{ i\in I:x_{i}\in
V\right\} $ is contained in $\mathcal{U}$ for any open neighborhood $V$ of $%
x $. This $\mathcal{U-}$limit is simply denoted by $\underset{i,\mathcal{U}}{%
\lim }x_{i}$. If a Banach space $E$ is infinite dimensional and $\mathcal{U}$
is an ultrafilter on the set $I$, then it is possible to enlarge $E$ to a
Banach space $\widehat{E}$ by adjoining for every bounded family $%
(x_{i})_{i\in I}$ in $E$ an element $\widehat{x}\in \widehat{E}$ such that
the following equality holds: $\left\Vert \widehat{x}\right\Vert =\underset{%
i,\mathcal{U}}{\lim }\left\Vert x_{i}\right\Vert $. This construction is
called \textit{Banach space ultrapower}. From the model-theoretical point of
view Banach space ultrapowers can be considered as the result of erazing the
elements with infinite norm from an ordinary $($i.e., algebraic$)$
ultrapower and dividing it by infinitesimals. Suppose that $(E_{i})_{i\in I}$
is a collection of Banach spaces. Then define%
\begin{equation*}
\ell _{\infty }\left( E_{i}\right) =\left\{ \left( x_{i}\right) :x_{i}\in
E_{i}\text{ and }\left\Vert \left( x_{i}\right) \right\Vert _{\infty
}<\infty \right\} \text{,}
\end{equation*}%
i.e., the Banach space of all bounded families $(x_{i})\in \underset{i\in I}{%
\dprod }E_{i}$ endowed with the norm given by $\left\Vert (x_{i})\right\Vert
_{\infty }=\underset{i\in I}{\sup }\left\Vert x_{i}\right\Vert _{E_{i}}$. If 
$\mathcal{U}$ is an ultrafilter on the index set $I$, then $\underset{i,%
\mathcal{U}}{\lim }\left\Vert x_{i}\right\Vert _{E_{i}}$ always exists. Thus
it is possible to define a seminorm on $\ell _{\infty }(E_{i})$ by assuming
that $\mathcal{N}((x_{i}))=\underset{i,\mathcal{U}}{\lim }\left\Vert
x_{i}\right\Vert _{E_{i}}$. Then the kernel of $\mathcal{N}$ is given by 
\begin{equation*}
\mathcal{N}_{\mathcal{U}}=\left\{ x=(x_{i})\in \ell _{\infty }\left(
E_{i}\right) :\underset{i,\mathcal{U}}{\lim }\left\Vert x_{i}\right\Vert
=0\right\}
\end{equation*}%
and it is straightforward to observe that $\mathcal{N}_{\mathcal{U}}$
constitutes a closed ideal in the Banach space $\ell _{\infty }\left(
E_{i}\right) $. Next, define the quotient space of the form%
\begin{equation*}
\ell _{\infty }\left( E_{i}\right) /\mathcal{N}_{\mathcal{U}}\text{.}
\end{equation*}%
This space is said to be the \textit{ultraproduct} of the family of Banach
spaces $\left( E_{i}\right) _{i\in I}$. If $E_{i}=E$ for every $i\in I$,
then the space $\ell _{\infty }(E_{i})/\mathcal{U}$ is termed the \textit{%
ultrapower} of $E$ and is denoted by $E^{I}/\mathcal{U}$ or by $\ell
_{\infty }(E)/\mathcal{N}_{\mathcal{U}}$. If $(x_{i})$ is the sequence in
the Banach space $\underset{i\in I}{\dprod }\left( E_{i}\right) $, then $%
(x_{i})_{\mathcal{U}}$ denotes the \textit{equivalence class} of $(x_{i})$
in $\ell _{\infty }(E_{i})/\mathcal{U}$. On the other hand, if $E^{I}/%
\mathcal{U}$ is an ultrapower of the Banach space $E$, then the mapping $%
x\rightarrow (x_{i})_{\mathcal{U}}$, where $x_{i}=x$ for every $i\in I$,
constitutes a classical embedding of $E$ into $E^{I}/\mathcal{U}$. It can be
asserted that if $\ell _{\infty }(E)/\mathcal{N}_{\mathcal{U}}$ is a Banach
space ultrapower, then $\ell _{\infty }(E)/\mathcal{N}_{\mathcal{U}}$
contains a subspace isometrically isomorphic to $E$ $($\cite{AK, AK1, BBHU,
BDCK, DCK, H, HI}$)$. \ \ \ \ \ \ \ \ \ \ \ \ \ \ \ \ \ \ \ \ \ \ \ \ \ \ \
\ \ \ \ \ \ \ \ \ \ \ \ \ \ \ \ \ \ \ \ \ \ \ \ \ \ \ \ \ \ \ \ \ \ \ \ \ \
\ \ \ \ \ \ \ \ \ \ \ \ \ \ \ \ \ \ \ \ \ \ \ \ \ 

A Banach space $E$ is said to be superreflexive if and only if each
ultrapower $\ell _{\infty }(E)/\mathcal{N}_{\mathcal{U}}$ is reflexive $($%
\cite{AK, H}$)$. There exist many alternative characterizations of
superreflexivity $($Theorem $2.1$ and $2.2$ in \cite{Wis}, cf. \cite{J}$)$.
In our paper the condition of superreflexivity of the Banach space $E$ will
be characterized by the containment of the copy of $\ell _{\infty }$ in $%
\ell _{\infty }(E)/\mathcal{N}_{\mathcal{U}}$.\ 

In $($\cite{Wil}$)$ the \textit{Representation Theorem} for nonreflexive
Banach space ultrapowers was obtained. In order to render to this result
recall that the \textit{Stone-Cech compactification} of the discrete space $%
\omega =\{1,2,...\}$, denoted by $\beta \omega $, can be identified with the
space of all ultrafilters on $\omega $ equipped with the topology generated
by sets of the form $\left\{ F:U\in F\right\} $ for $U\subseteq \omega $.
Then the points of $\beta \omega $ may be viewed as the ultrafilters on $%
\omega $ with the points of $\omega $ identified with the principal
ultrafilters. Denote by $\omega ^{\ast }$ $($or by $\beta \omega \backslash
\omega )$ the so-called \textit{Stone-Cech remainder }of $\beta \omega $. In
this case the points of $\omega ^{\ast }$ can be identified with free
ultrafilters on $\omega $ $($\cite{C, W}$)$. From \textit{Set-Theoretical
Topology} it is known that the so-called Parovicenco space $X$ can be
regarded as a topological space satisfying the following requirements: $1)$ $%
X$ is compact and Hausdorff, $2)$ $X$ has no isolated points, $3)$ $X$ has
weight $\mathfrak{c}$, $4)$ every nonempty $G_{\delta }$ subset of $X$ has
nonempty interior, $5)$ every two disjoint open $F_{\sigma }$ subsets of $X$
have disjoint closures. Parovi\v{c}enko proved that assuming $CH$ evey
Parovicenko space $X$ is isomorphic to $\omega ^{\ast }$ $($\cite{W}$)$. The
above mentioned Representation Theorem ascertains that - under $CH$ - all
nonreflexive Banach space ultrapowers can be represented in the form of the
space of continuous, bounded and real-valued functionsons the Parovicenco
space $\omega ^{\ast }$. Suppose that the symbol $\cong $ denotes the
relation of isometric isomorphism between Banach spaces. Then it is possible
to formulate the following theorem $($Corollary $3$ in \cite{Wil}$)$:\medskip

\textbf{Theorem }$\mathbf{1}$ $(CH)$. \textit{Let }$E$\textit{\ be any
infinite dimensional nonsuperreflexive Banach space and let }$\ell _{\infty
}(E)/N_{\mathcal{U}}$\textit{\ be its ultrapower. Then the Banach space }$%
\ell _{\infty }(E)/N_{\mathcal{U}}$\textit{\ is isometrically isomorphic to
the space of continuous, bounded and real-valued functions of the Stone-Cech
remainder }$\omega ^{\ast }$\textit{, i.e.,}%
\begin{equation*}
\ell _{\infty }(E)/\mathcal{N}_{\mathcal{U}}\cong C(\omega ^{\ast })\text{.}
\end{equation*}%
\medskip\ \ \ \ \ \ \ \ \ \ \ \ \ \ \ \ \ \ \ \ \ \ \ \ \ \ \ \ \ \ \ \ \ \
\ \ \ \ \ \ \ \ \ \ \ \ \ \ \ \ \ \ \ \ \ \ \ \ \ \ \ \ \ \ \ \ \ \ \ \ \ \
\ \ \ \ \ \ \ \ \ \ \ \ \ \ \ \ \ \ \ \ \ \ \ \ \ \ \ \ \ \ \ \ \ \ \ \ \ \
\ \ \ \ \ \ \ \ \ \ \ \ \ \ \ \ \ \ \ \ \ \ \ \ \ \ \ \ \ \ \ \ \ \ \ \ \ \
\ \ \ \ \ \ \ \ \ \ \ \ \ \ \ \ \ \ \ \ \ \ \ \ \ \ \ \ \ \ \ \ \ \ \ \ \ \
\ \ \ \ \ \ \ \ \ \ \ \ \ \ \ \ \ \ \ \ \ \ \ \ \ \ \ \ \ \ \ \ \ \ \ \ \ \
\ \ \ \ \ \ \ \ \ \ \ \ \ \ \ \ \ \ \ \ \ \ \ \ \ \ \ \ \ \ \ \ \ \ \ \ \ \
\ \ \ \ \ \ \ \ \ \ \ \ \ \ \ \ \ \ \ \ \ \ \ \ \ \ \ \ \ \ \ \ \ \ \ \ \ \
\ \ \ \ \ \ \ 

In order to outline the proof of this theorem suppose that the set $\omega $
has the discrete topology and $E$ is any infinite dimensional
nonsuperreflexive Banach space. Then $\ell _{\infty }(E)=C(\omega )$.
Defining the restriction mapping $R:C(\beta \omega )\rightarrow C(\omega )$
by $R(\widehat{f})=\widehat{f}\upharpoonright \omega $ for each $\widehat{f}%
\in C(\beta \omega )$ it can be concluded that $R$ is a linear isometry of $%
C(\beta \omega )$ onto $C(\omega )$. Consequently, it can be stated that if
we assume $CH$, then the $\ell _{\infty }-$sum of countably many copies of
any infinite dimensional nonsuperreflexive Banach space $E$ is isometrically
isomorphic to the space $C(\beta \omega )$, i.e., $\ell _{\infty }(E)\cong
C(\beta \omega )$ $($Proposition $1$ in \cite{Wil}$).$ Next, suppose that
every function $f\in C(\omega )$ has a unique norm-preserving extension $%
\widehat{f}\in C(\beta \omega )$. Therefore, it is possible to define the
closed ideal $I$ in the space $C(\beta \omega )$ consisting of functions
which vanish on $\omega ^{\ast }$, i.e., $I=\left\{ \widehat{f}\in C(\beta
\omega ):\widehat{f}(t)=0\text{ for all }t\in \omega ^{\ast }\right\} $.
Further, assume that the space $c_{0}(\omega )$ consists of functions in $%
C(\omega )$ which vanish at infinity, i.e., $c_{0}(\omega )=\left\{ f\in
C(\omega ):\text{ for each }\varepsilon >0\text{, }\left\{ t\in \omega
:\left\vert f(t)\right\vert >\varepsilon \right\} \text{ is finite}\right\} $%
. Observe that for $\underset{i,\mathcal{U}}{\lim }x_{i}=x$ the set $\left\{
i\in I:x_{i}\notin V\right\} $ is finite. Then it becomes obvious that $%
c_{0}(\omega )=\mathcal{N}_{\mathcal{U}}$ and the restriction mapping $%
R:I\rightarrow \mathcal{N}_{\mathcal{U}}$ defines a linear isometry from $I$
onto $\mathcal{N}_{\mathcal{U}}$. Consequently, $\mathcal{N}_{\mathcal{U}%
}\cong I$ $($Proposition $2$ in \cite{Wil}$)$. Summing up these two facts it
follows that $\ell _{\infty }(E)/\mathcal{N}_{\mathcal{U}}\cong C(\omega
^{\ast })$. \ \ \ \ \ \ \ \ \ \ \ \ \ \ \ \ \ \ \ \ \ \ \ \ \ \ \ \ \ \ \ \
\ \ \ \ \ \ \ \ \ \ \ \ \ \ \ \ \ \ \ \ \ \ \ \ \ \ \ \ \ \ \ \ \ \ \ \ \ \
\ \ \ \ \ \ \ \ \ \ \ \ \ \ \ \ \ \ \ \ \ \ \ \ \ \ \ \ \ \ \ \ \ \ \ \ \ \
\ \ \ \ \ \ \ \ \ \ \ \ \ \ \ \ \ \ \ \ \ \ \ \ \ \ \ \ \ \ \ \ \ \ \ \ \ \
\ \ \ \ \ \ \ \ \ \ \ \ \ \ \ \ \ \ \ \ \ \ \ \ \ \ \ \ \ \ \ \ \ \ \ \ \ \
\ \ \ \ \ \ \ \ \ \ \ \ \ \ \ \ \ \ \ \ \ \ \ \ \ \ \ \ \ \ \ \ \ \ \ \ \ \
\ \ \ \ \ \ \ \ \ \ \ \ \ \ \ \ \ \ \ \ \ \ \ \ \ \ \ \ \ \ \ \ \ \ \ \ \ \
\ \ \ \ \ \ \ \ \ \ \ \ \ \ \ \ \ \ \ \ \ \ \ \ \ \ \ \ \ \ \ \ \ \ \ \ \ \
\ \ \ 

In this place it shuold be noted that in our Representation Theorem it is
supposed that all considered infinite dimensional Banach spaces are
nonsuperreflexive and - consequently - their ultrapowers are nonreflexive.
It is unknown if the condition of nonsuperreflexivity can be weakened $($or
modified$)$ in order to represent Banach space ultrapowers in the form $%
C(\omega ^{\ast })$.

In this paper a large class of primary Banach spaces will be singled out.
Also it will be shown that - assuming $CH$ - every infinite dimensional
nonsuperreflexive Banach space can be isometrically embedded into its
primary ultrapower.\bigskip

$\mathbf{3.}$\textbf{\ All nonreflexive Banach space ultrapowers are primary
\bigskip }

Recall that if $X$ is any Banach space, then a sequence $\left\{
X_{n}\right\} _{n=1}^{\infty }$ of its closed subspaces is said to be a 
\textit{Schauder decomposition} of $X$ if every element $x\in X$ can be
uniquely represented in the form $x=\overset{\infty }{\underset{n=1}{\dsum }}
$, where $x_{n}\in X_{n}$ for every $n$. Drewnowski and Roberts showed that
if $CH$ holds, then there exists a Schauder decomposition of the space $%
C(\omega ^{\ast })$ $($Corollary $2.5$ in \cite{DR}, cf. \cite{D}$)$. Using
their result and the fact that the spaces $C(\omega ^{\ast })$ and $\ell
_{\infty }(E)/\mathcal{N}_{\mathcal{U}}$ $($where $E$ is any infinite
dimensional nonsuperreflexive Banach space$)$ are congruent the following
theorem can be formulated:\medskip

\textbf{Theorem }$\mathbf{2}$\textbf{.} \textit{Let }$E$\textit{\ be any
infinite dimensional nonsuperreflexive Banach space and let }$\ell _{\infty
}(E)/N_{\mathcal{U}}$\textit{\ be its ultrapower. Then for every }$($\textit{%
finite or infinite}$)$\textit{\ Schauder decomposition}%
\begin{equation*}
\ell _{\infty }(E)/\mathcal{N}_{\mathcal{U}}=X_{1}+X_{2}+...
\end{equation*}%
\textit{of Banach space ultrapower }$\ell _{\infty }(E)/N_{\mathcal{U}}$%
\textit{\ at least one of the summands }$X_{n}$\textit{\ has a subspace
which is isomorphic to }$\ell _{\infty }(E)/N_{\mathcal{U}}$\textit{\ and is
complemented in }$\ell _{\infty }(E)/N_{\mathcal{U}}$\textit{.}

\textit{Proof}. cf. \cite{DR}.\ $\square \medskip $

In order to completely grasp the primariness of Banach space ultrapowers let
us recall some fact from \textit{Set-Theoretical Topology}. Namely,
Negrepontis \ proved under $CH$ the following result $($Corollary $3.2$ in 
\cite{N}$)$:\medskip

\textbf{Proposition} $\mathbf{3}$ $(CH)$. \textit{If }$A$\textit{\ is an
open }$F_{\sigma }$\textit{\ subset of the Stone-Cech remainder }$\omega
^{\ast }$\textit{, then }$\overline{A}$\textit{\ constitutes a retract of }$%
\omega ^{\ast }.\medskip $

Now, suppose that the space $\ell _{\infty }\left( \ell _{\infty }\left(
E\right) /\mathcal{N}_{\mathcal{U}}\right) $ denotes the $\ell _{\infty }-$%
sum of countable many copies of Banach space ultrapower $\ell _{\infty
}\left( E\right) /\mathcal{N}_{\mathcal{U}}$, i.e., 
\begin{equation*}
\ell _{\infty }\left( \ell _{\infty }\left( E\right) /\mathcal{N}_{\mathcal{U%
}}\right) \colon =\left( \ell _{\infty }\left( E\right) /\mathcal{N}_{%
\mathcal{U}}\oplus \ell _{\infty }\left( E\right) /\mathcal{N}_{\mathcal{U}%
}\oplus ...\right) \text{.}
\end{equation*}%
In our next theorems we are going to show that if $E$ is any infinite
dimensional nonsuperreflexive Banach space, then the space $\ell _{\infty
}(\ell _{\infty }(E)/\mathcal{N}_{\mathcal{U}})$ is isomorphic to a
complemented subspace of Banach space ultrapower $\ell _{\infty }\left(
E\right) /\mathcal{N}_{\mathcal{U}}$. Consequently, it will be
straightforward to conclude that for any infinite dimensional
nonsuperreflexive Banach space $E$ the spaces $\ell _{\infty }(\ell _{\infty
}(E)/\mathcal{N}_{\mathcal{U}})$ and $\ell _{\infty }(E)/\mathcal{N}_{%
\mathcal{U}}$ are isomorphic, i.e., $\ell _{\infty }(\ell _{\infty }(E)/%
\mathcal{N}_{\mathcal{U}})\approx \ell _{\infty }(E)/\mathcal{N}_{\mathcal{U}%
}$. Also it should be noted that if $E$ is any infinite dimensional
nonsuperreflexive Banach space, then every infinite-dimensional complemented
subspace of $\ell _{\infty }(E)/\mathcal{N}_{\mathcal{U}}$ contains an
isomorph of the space $\ell _{\infty }(E)$ $($cf. \cite{D, DR}$)$.\medskip

\textbf{Theorem }$\mathbf{4}$ $(CH)$. \textit{Let }$E$\textit{\ be any
infinite dimensional nonsuperreflexive Banach space. Then the Banach space} $%
\ell _{\infty }(\ell _{\infty }(E)/\mathcal{N}_{\mathcal{U}})$ \textit{is
isometric to a complemented subspace of the Banach space ultrapower} $\ell
_{\infty }(E)/\mathcal{N}_{\mathcal{U}}$.

\textit{Proof. }Suppose that $\left\{ A_{n}\right\} _{n=1}^{\infty }$ is any
infinite sequence of disjoint and nonempty clopen subsets of the Parovicenco
space $\omega ^{\ast }$. Suppose that $A$ denotes the union of this
sequence. Then the Banach space $\ell _{\infty }(\ell _{\infty }(E)/\mathcal{%
N}_{\mathcal{U}})$ is isomorphically isometric to the Banach space $C_{b}(A)$
of bounded continuous functions on $A$. This fact follows easily from the
congruence between the spaces $\ell _{\infty }(E)/\mathcal{N}_{\mathcal{U}}$
and $C(\omega ^{\ast })$. Also it can be observed that $\overline{A}=\beta A$
and - consequently - $C_{b}\left( \overline{A}\right) \cong C(\overline{A})$%
. As the result it is obtained that $\ell _{\infty }(\ell _{\infty }(E)/%
\mathcal{N}_{\mathcal{U}})\cong C(\overline{A})$. From Proposition $3$ it is
deducible that there exists a retraction $r$ of $\omega ^{\ast }$ onto $%
\overline{A}$. If we define the corresponding composition operator $R:C(%
\overline{A})\rightarrow \ell _{\infty }(E)/\mathcal{N}_{\mathcal{U}}$, then
the mapping $x\rightarrow x\circ r$ can be identified with an isometry and
the operator $P$ given by $P(x)=R(x_{\mid A})$ is a projection \ from $\ell
_{\infty }(E)/\mathcal{N}_{\mathcal{U}}$ onto $R[C(\overline{A})]$. $\square
\medskip $

The proof of Theorem $5$ relies mainly on Pe\l czy\'{n}ski's decomposition
method $($cf. \cite{AK1, C, DR, LT}$)$.\medskip

\textbf{Theorem }$\mathbf{5}$\textbf{.} \textit{Let }$E$\textit{\ be any
infinite dimensional nonsuperreflexive Banach space. Then the Banach space }$%
\ell _{\infty }(\ell _{\infty }(E)/\mathcal{N}_{\mathcal{U}})$ \textit{is
isomorphic to the Banach space ultrapower }$\ell _{\infty }(E)/\mathcal{N}_{%
\mathcal{U}}$\textit{, i.e.,} 
\begin{equation*}
\ell _{\infty }(\ell _{\infty }(E)/\mathcal{N}_{\mathcal{U}})\approx \ell
_{\infty }(E)/\mathcal{N}_{\mathcal{U}}\text{.}
\end{equation*}%
\textit{Proof. }From Theorem\textit{\ }$4$ and Pe\l czy\'{n}ski's
decomposition technique it follows that if $\ell _{\infty }(E)/\mathcal{N}_{%
\mathcal{U}}\approx \ell _{\infty }(\ell _{\infty }(E)/\mathcal{N}_{\mathcal{%
U}})\oplus Z$ for some Banach space $Z$, then%
\begin{eqnarray*}
\ell _{\infty }(\ell _{\infty }(E)/\mathcal{N}_{\mathcal{U}}) &\approx &\ell
_{\infty }(\ell _{\infty }(E)/\mathcal{N}_{\mathcal{U}})\oplus \ell _{\infty
}(E)/\mathcal{N}_{\mathcal{U}} \\
&\approx &\ell _{\infty }(\ell _{\infty }(E)/\mathcal{N}_{\mathcal{U}%
})\oplus \ell _{\infty }(\ell _{\infty }(E)/\mathcal{N}_{\mathcal{U}})\oplus
Z \\
&\approx &\ell _{\infty }(\ell _{\infty }(E)/\mathcal{N}_{\mathcal{U}%
})\oplus Z \\
&\approx &\ell _{\infty }(E)/\mathcal{N}_{\mathcal{U}}\text{. }\square
\end{eqnarray*}%
\medskip From this result it is possible to deduce the additional
characterization of nonsuperreflexivity $($\cite{MA}$)$. Namely, the
following proposition can be formulated:\medskip

\textbf{Proposition }$\mathbf{6}$\textbf{.} \textit{If the Banach space }$%
\ell _{\infty }$\textit{\ of all bounded sequences is contained in the
ultrapower }$\ell _{\infty }(E)/N_{\mathcal{U}}$\textit{\ of an infinite
dimensional Banach space }$E$\textit{, then the Banach space }$E$\textit{\
is not superreflexive.\medskip }

Now we are ready to prove our main Theorem\ asserting that under $CH$ all
nonreflexive Banach space ultrapowers are primary $($cf. \cite{D, DR}$)$%
.\medskip

\textbf{Theorem }$\mathbf{7}$ $(CH).$ \textit{Let }$E$\textit{\ be any
infinite dimensional nonsuperreflexive Banach space and let }$\ell _{\infty
}(E)/\mathcal{N}_{\mathcal{U}}=X_{1}\oplus X_{2}\oplus ...$ \textit{be a
(finite or infinite) Schauder decomposition of its ultrapower. Then there
exists an index }$m$\textit{\ such that }$X_{m}\approx \ell _{\infty }(E)/%
\mathcal{N}_{\mathcal{U}}$.\textit{\ Particularly, the Banach space
ultrapower} $\ell _{\infty }(E)/\mathcal{N}_{\mathcal{U}}$ \textit{is
primary.}

\textit{Proof.} From Theorem $1$ it follows that it is possible to indicate
an index $m$ such that the space $X_{m}$ contains a subspace $V$ which is
isomorphic to $\ell _{\infty }(E)/\mathcal{N}_{\mathcal{U}}$ and is
complemented in $\ell _{\infty }(E)/\mathcal{N}_{\mathcal{U}}$. For
instance, suppose that $m=1$ and denote $X=X_{1}$, $Y=X_{2}\oplus
X_{3}\oplus ...$ . Hence we obtain:%
\begin{equation*}
\ell _{\infty }(E)/\mathcal{N}_{\mathcal{U}}=X\oplus Y\text{, }X=U\oplus V%
\text{ and }V\approx \ell _{\infty }(E)/\mathcal{N}_{\mathcal{U}}
\end{equation*}%
for some subspace $U$ of $X$. From Theorem $5$ we have that $\ell _{\infty
}(E)/\mathcal{N}_{\mathcal{U}}\approx \ell _{\infty }(\ell _{\infty }(E)/%
\mathcal{N}_{\mathcal{U}})$. Applying Pelczynski's decomposition method we
arrive at the following formula:%
\begin{eqnarray*}
X &\approx &U\oplus \ell _{\infty }(\ell _{\infty }(E)/\mathcal{N}_{\mathcal{%
U}})\approx U\oplus \ell _{\infty }(E)/\mathcal{N}_{\mathcal{U}}\oplus \ell
_{\infty }(\ell _{\infty }(E)/\mathcal{N}_{\mathcal{U}}) \\
&\approx &X\oplus \ell _{\infty }(X\oplus Y)\approx X\oplus \ell _{\infty
}(X)\oplus \ell _{\infty }(Y)\approx \ell _{\infty }(X)\oplus \ell _{\infty
}(Y) \\
&\approx &\ell _{\infty }(X\oplus Y)\approx \ell _{\infty }(\ell _{\infty
}(E)/\mathcal{N}_{\mathcal{U}})\approx \ell _{\infty }(E)/\mathcal{N}_{%
\mathcal{U}}\text{. }\square
\end{eqnarray*}%
$\bigskip \mathbf{4.}$\textbf{\ Concluding remarks}

In this paper it was demonstrated that if the \textit{Continuum Hypothesis}
holds, then all nonreflexive Banach space ultrapowers are primary. As the
consequence of this fact it can be observed that every infinite dimensional
nonsuperreflexive Banach space $E$ is isometrically embeddable into its
primary ultrapower $\ell _{\infty }(E)/\mathcal{N}_{\mathcal{U}}$. It was
also indicated that the presence of the copy of $\ell _{\infty }$ in $\ell
_{\infty }(E)/\mathcal{N}_{\mathcal{U}}$ $($where $E$ is infinite
dimensional Banach space$)$ implies that $E$ is not superreflexive.

\end{document}